\documentclass[12pt]{amsart}

\usepackage{fullpage}

\usepackage{pgf}
\usepackage{graphicx}
\usepackage{amsmath, amsthm, amssymb, amscd}
\usepackage{hyperref}

\usepackage{tikz}

\def\hrad{0.5}
\def\vrad{0.1}
\def\rsep{0.25}
\def\miny{0}
\def\maxy{4}
\def\tension{0.5}

\newcommand\Z{{\mathbb Z}}

\newcommand\R{{\mathbb R}}
\newcommand\C{{\mathbb C}}

\newcommand\Q{{\mathbb Q}}

\newcommand\bp{\begin{proof}}
\newcommand\ep{\end{proof}}

\newcommand\bprop{\begin{prop}}
\newcommand\eprop{\end{prop}}

\newcommand\bt{\begin{thm}}
\newcommand\et{\end{thm}}

\newcommand\bc{\begin{cor}}
\newcommand\ec{\end{cor}}

\newcommand\ba{\begin{aligned}}
\newcommand\ea{\end{aligned}}

\newcommand\bl{\begin{lem}}
\newcommand\el{\end{lem}}

\newcommand\bi{\begin{itemize}}
\newcommand\ei{\end{itemize}}

\newcommand\br{\begin{rem}}
\newcommand\er{\end{rem}}

\newcommand\bd{\begin{defn}}
\newcommand\ed{\end{defn}}

\newcommand\B{{\mathcal B}}

\newcommand{\mcLB}{\mathcal{LB} }
\newcommand{\LB}{\mathcal{LB} }
\newcommand{\SLB}{\mathcal{SLB} }
\newcommand{\VB}{\mathcal{VB} }
\newcommand{\ot}{\otimes }

\DeclareMathOperator{\Hom}{Hom}

\DeclareMathOperator{\End}{End}

\DeclareMathOperator{\diag}{diag}

\DeclareMathOperator{\ind}{ind}
\DeclareMathOperator{\res}{res}

\DeclareMathOperator{\Gr}{Gr}
\DeclareMathOperator{\Det}{Det}

\DeclareMathOperator{\GL}{GL}

\DeclareMathOperator{\U}{U}
\newtheorem{thm}{Theorem}[section]

\newtheorem*{thm*}{Theorem}

\newtheorem{prop}[thm]{Proposition}
\newtheorem{lem}[thm]{Lemma}
\newtheorem{cor}[thm]{Corollary}

\newtheorem*{q*}{Question}

\theoremstyle{remark}
\newtheorem{rem}[thm]{Remark}
\newtheorem*{rem*}{Remark}

\theoremstyle{definition}
\newtheorem{defn}[thm]{Definition}

\newtheorem{conj}{Conjecture}
\newtheorem{eg}[thm]{Example}

\renewcommand{\l}{\lambda}

\newcommand{\Tr}{\mathrm{Tr}}

\author{Paul Bruillard}
\email{pjb2357@gmail.com}
\address{Pacific Northwest National Laboratory, 902 Battelle Boulevard,
Richland, WA U.S.A}
\thanks{\textit{PNNL Information Release:} PNNL-SA-111814.}

    \author{Seung-Moon Hong}
\email{seungmoon.hong@utoledo.edu}
\address{Department of Mathematics and Statistics,
    University of Toledo,
    Toledo, OH
    U.S.A.}
\author{Julia Yael Plavnik}
\email{plavnik@famaf.unc.edu.ar}
\address{Facultad de Matem\'atica, Astronom\'ia y F\'isica,
Universidad Nacional de C\'ordoba, CIEM -- CONICET, (5000)
 Ciudad Universitaria, C\'ordoba, Argentina.}

\author{Eric C. Rowell}
\email{rowell@math.tamu.edu}
\address{Department of Mathematics,
    Texas A\&M University,
    College Station, TX
    U.S.A.}
\author{Liang Chang}
    \email{changliang996@gmail.com}
    \address{Department of Mathematics,
    Texas A\&M University,
    College Station, TX
    U.S.A.}

\author{Michael Yuan Sun}
\email{sunm@uni-muenster.de}
\address{Mathematisches Institut der Universit\"{a}t M\"{u}nster, Germany}

\thanks{This paper began while the authors were participating in an AMS
Mathematics Research Community workshop in Snowbird, Utah.  The continued
support of the AMS is gratefully acknowledged. J.Y.P was partially supported by
CONICET, ANPCyT  and Secyt (UNC). E.C.R. was partially supported by NSF
DMS-1108725. M.Y.S. acknowledges the support of SFB878 and GIF.
The research described in this paper was, in part, conducted under the
Laboratory Directed Research and Development Program at PNNL, a multi-program
national laboratory operated by Battelle for the U.S. Department of Energy.
}
\title[Low-dimensional representations of $\mathcal{LB}_3$]{Low-dimensional representations of the three component loop braid group}
\begin{document}
\begin{abstract}

Motivated by physical and topological applications, we study representations of the group $\LB_3$ of motions of $3$ unlinked oriented circles in $\R^3$.  Our point of view is to regard the three strand braid group $\B_3$ as a subgroup of $\LB_3$ and study the problem of extending $\B_3$ representations.  We introduce the notion of a \emph{standard extension} and characterize $\B_3$ representations admiting such an extension.  In particular we show, using a classification result of Tuba and Wenzl, that every irreducible $\B_3$ representation of dimension at most $5$ has a (standard) extension.
We show that this result is sharp by exhibiting an irreducible $6$-dimensional $\B_3$ representation that has no extensions (standard or otherwise).  We obtain complete classifications of (1) irreducible $2$-dimensional $\LB_3$ representations (2) extensions of irreducible $3$-dimensional $\B_3$ representations and (3) irreducible $\LB_3$ representations whose restriction to $\B_3$ has abelian image.

\end{abstract}
\maketitle
\section{Introduction and Motivation}
Over the last two decades, topological states of matter in $2$ spatial dimensions and their potential computational applications have motivated the study of motions of sytems of point-like excitations on $2$-dimensional surfaces.  The mathematical model for such systems are $(2+1)$-topological quantum field theories (TQFTs).  The motions of points in the disk lead to representations of the braid group $\B_n$ which play a central role in the topological model for quantum computation \cite{FKLW,FLW1}.

Recently the possibility of $3$-dimensional topological states of matter (see, e.g. \cite{WL,WW}) modeled by $(3+1)$-TQFTs \cite{ww12} lead to new possible avenues for quantum computation.  Although motions of point-like excitations in $3$ spatial dimensions are mathematically trivial, the symmetries of loop-like excitations can be quite complicated.  The simplest mathematical manifestation of this idea is the group $\mcLB_n$ of motions of an oriented $n$-component unlink $C_n$ \cite{dahm,Gold,lin}.
In this article, we begin a systematic study of the low-dimensional representations of the loop braid group $\mcLB_3$.  The generators and relations for $\mcLB_3$ are given in \cite{FRR}. We take it as a definition here:
\begin{defn}
The three component \textbf{loop braid group} $\mcLB_3$ is the abstract group generated by $\sigma_1,\sigma_2,s_1,s_2$ satisfying the following relations:
\begin{enumerate}
 \item[(B1)] $\sigma_1\sigma_{2}\sigma_1=\sigma_{2}\sigma_1\sigma_{2}$
 \item[($\mathfrak{S}$1)] $s_1s_{2}s_1=s_{2}s_1s_{2}$
 \item[($\mathfrak{S}$2)] $s_1^2=s_2^2=1$
\item[(L1)] $s_1s_{2}\sigma_1=\sigma_{2}s_1s_{2}$
\item[(L2)] $\sigma_1\sigma_{2}s_1=s_{2}\sigma_1\sigma_{2}$
\end{enumerate}
\end{defn}
The group satisfying only (B1), ($\mathfrak{S}$1) and ($\mathfrak{S}$2) is the free product $\B_3\ast\mathfrak{S}_3$ of Artin's braid group $\B_3$ and the symmetric group $\mathfrak{S}_3$.  The group satisfying (B1)-(L1) is known as the virtual braid group on three strands $\mathcal{VB}_3$ \cite{Kauf}.
Some authors replace (L2) by
\begin{enumerate}
 \item[(L2$^\prime$)] $\sigma_2\sigma_{1}s_2=s_{1}\sigma_2\sigma_{1}$
\end{enumerate}
which leads to an isomorphic group. The group satisfying (B1)-(L2) \emph{and} (L2$^\prime$) is called the \textbf{symmetric loop braid group} $\mathcal{SLB}_3$ \cite{leeds}.

The $3$-component loop braid group $\mcLB_3$ is geometrically understood as motions of $3$ oriented circles in $\R^3$. The generator $\sigma_i$ is interpreted as passing the $i$th circle under and through the $i+1$st circle ending with the two circles' positions interchanged. The generator $s_i$ corresponds to simply interchanging circles $i$ and $i+1$.  We can represent these generators diagrammatically  as follows, with the time variable to be read from bottom to top:


  
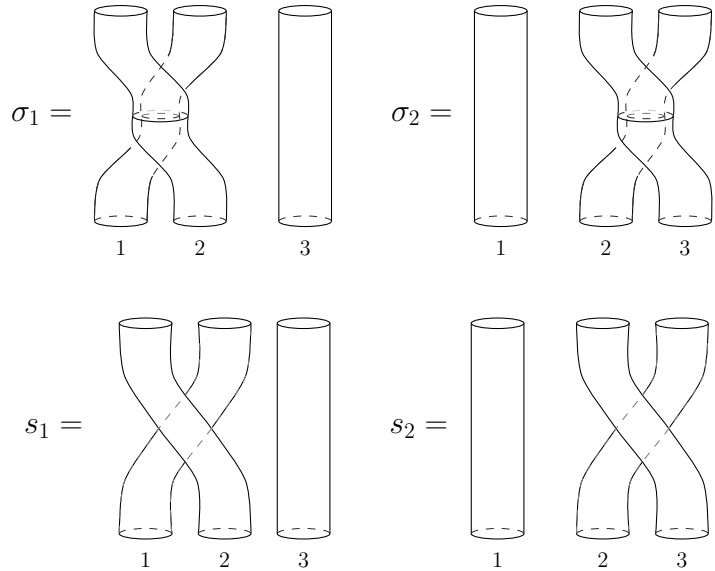
\begin{center}
 \begin{tikzpicture}[scale=0.7,every node/.style={scale=0.7}]
      \node[] at ({-4*\rsep},{\maxy/2}) {\Large{$\sigma_{1} = $}};
      \begin{scope}[xshift = 3.5cm]
        \draw (0,0) arc (180:360:{\hrad} and {\vrad});
        \draw[dashed,color=black!80!white] (0,0) arc (180:0:{\hrad} and {\vrad});
        \draw (0,{\miny}) -- (0,{\maxy});
        \draw (1,{\miny}) -- (1,{\maxy});
        \draw (0,{\maxy}) arc(180:360:{\hrad} and {\vrad});
        \draw (0,{\maxy}) arc(180:0:{\hrad} and {\vrad});
        \node[] at ({\hrad},{\miny-2*\rsep}) {$3$};
      \end{scope}
      \begin{scope}[xshift = 0cm]
        \draw plot [smooth,tension={\tension}] coordinates{
        ({\hrad-\hrad},0)
        ({\hrad-0.8*\hrad},{\maxy/5})
        ({2*\hrad+\rsep-0.8*\hrad},{2*\maxy/5})
        ({2*\hrad+\rsep-0.7*\hrad},{3*\maxy/5})
        ({3*\hrad+2*\rsep-1.15*\hrad},{4*\maxy/5})
        ({3*\hrad+2*\rsep-\hrad},{5*\maxy/5})};
        \draw plot [smooth,tension={\tension}] coordinates{
        ({\hrad+\hrad},0)
        ({\hrad+1.15*\hrad},{\maxy/5})
        ({2*\hrad+\rsep+0.7*\hrad},{2*\maxy/5})
        ({2*\hrad+\rsep+0.8*\hrad},{3*\maxy/5})
        ({3*\hrad+2*\rsep+0.8*\hrad},{4*\maxy/5})
        ({3*\hrad+2*\rsep+\hrad},{5*\maxy/5})};
        \draw[line width={60*0.8*\hrad},color=white] plot [smooth,tension={\tension}] coordinates{
        ({3*\hrad+2*\rsep},0)
        ({3*\hrad+2*\rsep},{\maxy/5})
        ({2*\hrad+\rsep},{2*\maxy/5})
        ({2*\hrad+\rsep},{3*\maxy/5})
        ({\hrad},{4*\maxy/5})
        ({\hrad},{5*\maxy/5})};
        \draw plot [smooth,tension={\tension}] coordinates{
        ({3*\hrad+2*\rsep-\hrad},0)
        ({3*\hrad+2*\rsep-1.1*\hrad},{\maxy/5})
        ({2*\hrad+\rsep-\hrad},{2*\maxy/5})
        ({2*\hrad+\rsep-1.1*\hrad},{3*\maxy/5})
        ({\hrad-0.9*\hrad},{4*\maxy/5})
        ({\hrad-\hrad},{5*\maxy/5})};
        \draw plot [smooth,tension={\tension}] coordinates{
        ({3*\hrad+2*\rsep+\hrad},0)
        ({3*\hrad+2*\rsep+0.9*\hrad},{\maxy/5})
        ({2*\hrad+\rsep+1.1*\hrad},{2*\maxy/5})
        ({2*\hrad+\rsep+\hrad},{3*\maxy/5})
        ({\hrad+1.1*\hrad},{4*\maxy/5})
        ({\hrad+\hrad},{5*\maxy/5})};
        \draw[dashed,color=black!80!white] plot [smooth,tension={\tension}] coordinates{
        ({\hrad-\hrad},0)
        ({\hrad-0.8*\hrad},{\maxy/5})
        ({2*\hrad+\rsep-0.8*\hrad},{2*\maxy/5})
        ({2*\hrad+\rsep-0.7*\hrad},{3*\maxy/5})
        ({3*\hrad+2*\rsep-1.15*\hrad},{4*\maxy/5})
        ({3*\hrad+2*\rsep-\hrad},{5*\maxy/5})};
        \draw[dashed, color=black!80!white] plot [smooth,tension={\tension}] coordinates{
        ({\hrad+\hrad},0)
        ({\hrad+1.15*\hrad},{\maxy/5})
        ({2*\hrad+\rsep+0.7*\hrad},{2*\maxy/5})
        ({2*\hrad+\rsep+0.8*\hrad},{3*\maxy/5})
        ({3*\hrad+2*\rsep+0.8*\hrad},{4*\maxy/5})
        ({3*\hrad+2*\rsep+\hrad},{5*\maxy/5})};
        \draw ({0.95*\hrad+\rsep},{2.5*\maxy/5}) arc (180:360:{1.05*\hrad} and {1.1*\vrad});
        \draw[dashed,color=black!20!white] ({0.95*\hrad+\rsep},{2.5*\maxy/5}) arc (180:0:{1.05*\hrad} and {1.1*\vrad});
        \draw[color=white] ({0.95*\hrad+\rsep},{2.5*\maxy/5}) arc (180:130:{1.05*\hrad} and {1.1*\vrad});
        \draw[dashed,color=black!80!white] ({0.95*\hrad+\rsep},{2.5*\maxy/5}) arc (180:130:{1.05*\hrad} and {1.1*\vrad});
        \draw[color=white] ({0.95*\hrad+\rsep+2.1*\hrad},{2.5*\maxy/5}) arc (0:50:{1.05*\hrad} and {1.1*\vrad});
        \draw[dashed,color=black!80!white] ({0.95*\hrad+\rsep+2.1*\hrad},{2.5*\maxy/5}) arc (0:50:{1.05*\hrad} and {1.1*\vrad});
        \draw[dashed,color=black!80!white] ({2*\hrad+\rsep-0.7*\hrad},{2.5*\maxy/5}) arc (180:360:{0.7*\hrad} and {0.5*\vrad});
        \draw[dashed,color=black!50!white] ({2*\hrad+\rsep-0.7*\hrad},{2.5*\maxy/5}) arc (180:0:{0.7*\hrad} and {0.5*\vrad});
        \draw (0,0) arc (180:360:{\hrad} and {\vrad});
        \draw[dashed,color=black!80!white] (0,0) arc (180:0:{\hrad} and {\vrad});
        \draw (0,{\maxy}) arc(180:360:{\hrad} and {\vrad});
        \draw (0,{\maxy}) arc(180:0:{\hrad} and {\vrad});
        \draw ({2*\hrad+2*\rsep},0) arc (180:360:{\hrad} and {\vrad});
        \draw[dashed,color=black!80!white] ({2*\hrad+2*\rsep},0) arc (180:0:{\hrad} and {\vrad});
        \draw ({2*\hrad+2*\rsep},{\maxy}) arc(180:360:{\hrad} and {\vrad});
        \draw ({2*\hrad+2*\rsep},{\maxy}) arc(180:0:{\hrad} and {\vrad});
        \node[] at ({\hrad},{\miny-2*\rsep}) {$1$};
        \node[] at ({2*\hrad+2*\rsep+\hrad},{\miny-2*\rsep}) {$2$};
      \end{scope}
    \end{tikzpicture}
    $\quad$
       \begin{tikzpicture}[scale=0.7,every node/.style={scale=0.7}]
      \node[] at ({-4*\rsep},{\maxy/2}) {\Large{$\sigma_{2} = $}};
      \begin{scope}[xshift = 0]
        \draw (0,0) arc (180:360:{\hrad} and {\vrad});
        \draw[dashed,color=black!80!white] (0,0) arc (180:0:{\hrad} and {\vrad});
        \draw (0,{\miny}) -- (0,{\maxy});
        \draw (1,{\miny}) -- (1,{\maxy});
        \draw (0,{\maxy}) arc(180:360:{\hrad} and {\vrad});
        \draw (0,{\maxy}) arc(180:0:{\hrad} and {\vrad});
        \node[] at ({\hrad},{\miny-2*\rsep}) {$1$};
      \end{scope}
      \begin{scope}[xshift = 2cm]
        \draw plot [smooth,tension={\tension}] coordinates{
        ({\hrad-\hrad},0)
        ({\hrad-0.8*\hrad},{\maxy/5})
        ({2*\hrad+\rsep-0.8*\hrad},{2*\maxy/5})
        ({2*\hrad+\rsep-0.7*\hrad},{3*\maxy/5})
        ({3*\hrad+2*\rsep-1.15*\hrad},{4*\maxy/5})
        ({3*\hrad+2*\rsep-\hrad},{5*\maxy/5})};
        \draw plot [smooth,tension={\tension}] coordinates{
        ({\hrad+\hrad},0)
        ({\hrad+1.15*\hrad},{\maxy/5})
        ({2*\hrad+\rsep+0.7*\hrad},{2*\maxy/5})
        ({2*\hrad+\rsep+0.8*\hrad},{3*\maxy/5})
        ({3*\hrad+2*\rsep+0.8*\hrad},{4*\maxy/5})
        ({3*\hrad+2*\rsep+\hrad},{5*\maxy/5})};
        \draw[line width={60*0.8*\hrad},color=white] plot [smooth,tension={\tension}] coordinates{
        ({3*\hrad+2*\rsep},0)
        ({3*\hrad+2*\rsep},{\maxy/5})
        ({2*\hrad+\rsep},{2*\maxy/5})
        ({2*\hrad+\rsep},{3*\maxy/5})
        ({\hrad},{4*\maxy/5})
        ({\hrad},{5*\maxy/5})};
        \draw plot [smooth,tension={\tension}] coordinates{
        ({3*\hrad+2*\rsep-\hrad},0)
        ({3*\hrad+2*\rsep-1.1*\hrad},{\maxy/5})
        ({2*\hrad+\rsep-\hrad},{2*\maxy/5})
        ({2*\hrad+\rsep-1.1*\hrad},{3*\maxy/5})
        ({\hrad-0.9*\hrad},{4*\maxy/5})
        ({\hrad-\hrad},{5*\maxy/5})};
        \draw plot [smooth,tension={\tension}] coordinates{
        ({3*\hrad+2*\rsep+\hrad},0)
        ({3*\hrad+2*\rsep+0.9*\hrad},{\maxy/5})
        ({2*\hrad+\rsep+1.1*\hrad},{2*\maxy/5})
        ({2*\hrad+\rsep+\hrad},{3*\maxy/5})
        ({\hrad+1.1*\hrad},{4*\maxy/5})
        ({\hrad+\hrad},{5*\maxy/5})};
        \draw[dashed,color=black!80!white] plot [smooth,tension={\tension}] coordinates{
        ({\hrad-\hrad},0)
        ({\hrad-0.8*\hrad},{\maxy/5})
        ({2*\hrad+\rsep-0.8*\hrad},{2*\maxy/5})
        ({2*\hrad+\rsep-0.7*\hrad},{3*\maxy/5})
        ({3*\hrad+2*\rsep-1.15*\hrad},{4*\maxy/5})
        ({3*\hrad+2*\rsep-\hrad},{5*\maxy/5})};
        \draw[dashed, color=black!80!white] plot [smooth,tension={\tension}] coordinates{
        ({\hrad+\hrad},0)
        ({\hrad+1.15*\hrad},{\maxy/5})
        ({2*\hrad+\rsep+0.7*\hrad},{2*\maxy/5})
        ({2*\hrad+\rsep+0.8*\hrad},{3*\maxy/5})
        ({3*\hrad+2*\rsep+0.8*\hrad},{4*\maxy/5})
        ({3*\hrad+2*\rsep+\hrad},{5*\maxy/5})};
        \draw ({0.95*\hrad+\rsep},{2.5*\maxy/5}) arc (180:360:{1.05*\hrad} and {1.1*\vrad});
        \draw[dashed,color=black!20!white] ({0.95*\hrad+\rsep},{2.5*\maxy/5}) arc (180:0:{1.05*\hrad} and {1.1*\vrad});
        \draw[color=white] ({0.95*\hrad+\rsep},{2.5*\maxy/5}) arc (180:130:{1.05*\hrad} and {1.1*\vrad});
        \draw[dashed,color=black!80!white] ({0.95*\hrad+\rsep},{2.5*\maxy/5}) arc (180:130:{1.05*\hrad} and {1.1*\vrad});
        \draw[color=white] ({0.95*\hrad+\rsep+2.1*\hrad},{2.5*\maxy/5}) arc (0:50:{1.05*\hrad} and {1.1*\vrad});
        \draw[dashed,color=black!80!white] ({0.95*\hrad+\rsep+2.1*\hrad},{2.5*\maxy/5}) arc (0:50:{1.05*\hrad} and {1.1*\vrad});
        \draw[dashed,color=black!80!white] ({2*\hrad+\rsep-0.7*\hrad},{2.5*\maxy/5}) arc (180:360:{0.7*\hrad} and {0.5*\vrad});
        \draw[dashed,color=black!50!white] ({2*\hrad+\rsep-0.7*\hrad},{2.5*\maxy/5}) arc (180:0:{0.7*\hrad} and {0.5*\vrad});
        \draw (0,0) arc (180:360:{\hrad} and {\vrad});
        \draw[dashed,color=black!80!white] (0,0) arc (180:0:{\hrad} and {\vrad});
        \draw (0,{\maxy}) arc(180:360:{\hrad} and {\vrad});
        \draw (0,{\maxy}) arc(180:0:{\hrad} and {\vrad});
        \draw ({2*\hrad+2*\rsep},0) arc (180:360:{\hrad} and {\vrad});
        \draw[dashed,color=black!80!white] ({2*\hrad+2*\rsep},0) arc (180:0:{\hrad} and {\vrad});
        \draw ({2*\hrad+2*\rsep},{\maxy}) arc(180:360:{\hrad} and {\vrad});
        \draw ({2*\hrad+2*\rsep},{\maxy}) arc(180:0:{\hrad} and {\vrad});
        \node[] at ({\hrad},{\miny-2*\rsep}) {$2$};
        \node[] at ({2*\hrad+2*\rsep+\hrad},{\miny-2*\rsep}) {$3$};
      \end{scope}
    \end{tikzpicture}

    \ \\

    \begin{tikzpicture}[scale = 0.7,every node/.style={scale=0.7}]
      \node[] at ({-3*\rsep},{\maxy/2}) {\Large{$s_{1} = $}};
      \begin{scope}[xshift = 3.5cm]
        \draw (0,0) arc (180:360:{\hrad} and {\vrad});
        \draw[dashed,color=black!80!white] (0,0) arc (180:0:{\hrad} and {\vrad});
        \draw (0,{\miny}) -- (0,{\maxy});
        \draw (1,{\miny}) -- (1,{\maxy});
        \draw (0,{\maxy}) arc(180:360:{\hrad} and {\vrad});
        \draw (0,{\maxy}) arc(180:0:{\hrad} and {\vrad});
        \node[] at ({\hrad},{\miny-2*\rsep}) {$3$};
      \end{scope}

      \begin{scope}[xshift = 3cm,xscale=-1]
        \draw plot [smooth,tension={\tension}] coordinates{
        ({3*\hrad+2*\rsep-\hrad},0)
        ({3*\hrad+2*\rsep-1.05*\hrad},{\maxy/4})
        ({2*\hrad+\rsep-\hrad},{\maxy/2})
        ({\hrad-0.8*\hrad},{3*\maxy/4})
        ({\hrad-\hrad},{\maxy})};
        \draw plot [smooth,tension={\tension}] coordinates{
        ({3*\hrad+2*\rsep+\hrad},0)
        ({3*\hrad+2*\rsep+0.8*\hrad},{\maxy/4})
        ({2*\hrad+\rsep+\hrad},{\maxy/2})
        ({\hrad+1.05*\hrad},{3*\maxy/4})
        ({\hrad+\hrad},{\maxy})};
        \draw[color=white,line width={60*0.55*\hrad}] plot [smooth,tension={\tension}] coordinates{({\hrad},0) ({\hrad},{\maxy/4}) ({2*\hrad+\rsep},{\maxy/2}) ({3*\hrad+2*\rsep},{3*\maxy/4}) ({3*\hrad+2*\rsep},{\maxy})};
        \draw plot [smooth,tension={\tension}] coordinates{
        ({\hrad-\hrad},0)
        ({\hrad-0.8*\hrad},{\maxy/4})
        ({2*\hrad+\rsep-\hrad},{\maxy/2})
        ({3*\hrad+2*\rsep-1.05*\hrad},{3*\maxy/4})
        ({3*\hrad+2*\rsep-\hrad},{\maxy})};
        \draw plot [smooth,tension={\tension}] coordinates{
        ({\hrad+\hrad},0)
        ({\hrad+1.05*\hrad},{\maxy/4})
        ({2*\hrad+\rsep+\hrad},{\maxy/2})
        ({3*\hrad+2*\rsep+0.8*\hrad},{3*\maxy/4})
        ({3*\hrad+2*\rsep+\hrad},{\maxy})};
        \draw[dashed,color=black!50!white] plot [smooth,tension={\tension}] coordinates{
        ({3*\hrad+2*\rsep-\hrad},0)
        ({3*\hrad+2*\rsep-1.05*\hrad},{\maxy/4})
        ({2*\hrad+\rsep-\hrad},{\maxy/2})
        ({\hrad-0.8*\hrad},{3*\maxy/4})
        ({\hrad-\hrad},{\maxy})};
        \draw[dashed,color=black!50!white] plot [smooth,tension={\tension}] coordinates{
        ({3*\hrad+2*\rsep+\hrad},0)
        ({3*\hrad+2*\rsep+0.8*\hrad},{\maxy/4})
        ({2*\hrad+\rsep+\hrad},{\maxy/2})
        ({\hrad+1.05*\hrad},{3*\maxy/4})
        ({\hrad+\hrad},{\maxy})};
        \draw (0,0) arc (180:360:{\hrad} and {\vrad});
        \draw[dashed,color=black!80!white] (0,0) arc (180:0:{\hrad} and {\vrad});
        \draw (0,{\maxy}) arc(180:360:{\hrad} and {\vrad});
        \draw (0,{\maxy}) arc(180:0:{\hrad} and {\vrad});
        \draw ({2*\hrad+2*\rsep},0) arc (180:360:{\hrad} and {\vrad});
        \draw[dashed,color=black!80!white] ({2*\hrad+2*\rsep},0) arc (180:0:{\hrad} and {\vrad});
        \draw ({2*\hrad+2*\rsep},{\maxy}) arc(180:360:{\hrad} and {\vrad});
        \draw ({2*\hrad+2*\rsep},{\maxy}) arc(180:0:{\hrad} and {\vrad});
      \end{scope}
        \node[] at ({2*\hrad},{\miny-2*\rsep}) {$1$};
        \node[] at ({3*\hrad+2*\rsep+\hrad},{\miny-2*\rsep}) {$2$};
    \end{tikzpicture}
     $\quad$
    \begin{tikzpicture}[scale = 0.7,every node/.style={scale=0.7}]
      \node[] at ({-4*\rsep},{\maxy/2}) {\Large{$s_{2} = $}};
      \begin{scope}[xshift = 0]
        \draw (0,0) arc (180:360:{\hrad} and {\vrad});
        \draw[dashed,color=black!80!white] (0,0) arc (180:0:{\hrad} and {\vrad});
        \draw (0,{\miny}) -- (0,{\maxy});
        \draw (1,{\miny}) -- (1,{\maxy});
        \draw (0,{\maxy}) arc(180:360:{\hrad} and {\vrad});
        \draw (0,{\maxy}) arc(180:0:{\hrad} and {\vrad});
        \node[] at ({\hrad},{\miny-2*\rsep}) {$1$};
      \end{scope}

      \begin{scope}[xshift = 4.5cm,xscale=-1]
        \draw plot [smooth,tension={\tension}] coordinates{
        ({3*\hrad+2*\rsep-\hrad},0)
        ({3*\hrad+2*\rsep-1.05*\hrad},{\maxy/4})
        ({2*\hrad+\rsep-\hrad},{\maxy/2})
        ({\hrad-0.8*\hrad},{3*\maxy/4})
        ({\hrad-\hrad},{\maxy})};
        \draw plot [smooth,tension={\tension}] coordinates{
        ({3*\hrad+2*\rsep+\hrad},0)
        ({3*\hrad+2*\rsep+0.8*\hrad},{\maxy/4})
        ({2*\hrad+\rsep+\hrad},{\maxy/2})
        ({\hrad+1.05*\hrad},{3*\maxy/4})
        ({\hrad+\hrad},{\maxy})};
        \draw[color=white,line width={60*0.55*\hrad}] plot [smooth,tension={\tension}] coordinates{({\hrad},0) ({\hrad},{\maxy/4}) ({2*\hrad+\rsep},{\maxy/2}) ({3*\hrad+2*\rsep},{3*\maxy/4}) ({3*\hrad+2*\rsep},{\maxy})};
        \draw plot [smooth,tension={\tension}] coordinates{
        ({\hrad-\hrad},0)
        ({\hrad-0.8*\hrad},{\maxy/4})
        ({2*\hrad+\rsep-\hrad},{\maxy/2})
        ({3*\hrad+2*\rsep-1.05*\hrad},{3*\maxy/4})
        ({3*\hrad+2*\rsep-\hrad},{\maxy})};
        \draw plot [smooth,tension={\tension}] coordinates{
        ({\hrad+\hrad},0)
        ({\hrad+1.05*\hrad},{\maxy/4})
        ({2*\hrad+\rsep+\hrad},{\maxy/2})
        ({3*\hrad+2*\rsep+0.8*\hrad},{3*\maxy/4})
        ({3*\hrad+2*\rsep+\hrad},{\maxy})};
        \draw[dashed,color=black!50!white] plot [smooth,tension={\tension}] coordinates{
        ({3*\hrad+2*\rsep-\hrad},0)
        ({3*\hrad+2*\rsep-1.05*\hrad},{\maxy/4})
        ({2*\hrad+\rsep-\hrad},{\maxy/2})
        ({\hrad-0.8*\hrad},{3*\maxy/4})
        ({\hrad-\hrad},{\maxy})};
        \draw[dashed,color=black!50!white] plot [smooth,tension={\tension}] coordinates{
        ({3*\hrad+2*\rsep+\hrad},0)
        ({3*\hrad+2*\rsep+0.8*\hrad},{\maxy/4})
        ({2*\hrad+\rsep+\hrad},{\maxy/2})
        ({\hrad+1.05*\hrad},{3*\maxy/4})
        ({\hrad+\hrad},{\maxy})};
        \draw (0,0) arc (180:360:{\hrad} and {\vrad});
        \draw[dashed,color=black!80!white] (0,0) arc (180:0:{\hrad} and {\vrad});
        \draw (0,{\maxy}) arc(180:360:{\hrad} and {\vrad});
        \draw (0,{\maxy}) arc(180:0:{\hrad} and {\vrad});
        \draw ({2*\hrad+2*\rsep},0) arc (180:360:{\hrad} and {\vrad});
        \draw[dashed,color=black!80!white] ({2*\hrad+2*\rsep},0) arc (180:0:{\hrad} and {\vrad});
        \draw ({2*\hrad+2*\rsep},{\maxy}) arc(180:360:{\hrad} and {\vrad});
        \draw ({2*\hrad+2*\rsep},{\maxy}) arc(180:0:{\hrad} and {\vrad});
      \end{scope}
        \node[] at ({2+\hrad},{\miny-2*\rsep}) {$2$};
        \node[] at ({2+2*\hrad+2*\rsep+\hrad},{\miny-2*\rsep}) {$3$};
    \end{tikzpicture}

  \end{center}

A few authors have approached the representations of $\mcLB_3$ from the point of view of extending representations of $\B_3$ (see \cite{BWC,bard1,bard2,leeds,Ver}).  Among these, the extending of specific families of representations of $\B_3$, such as the Burau and Lawrence-Krammer-Bigelow representations, as well as certain representations obtained from solutions to the Yang-Baxter equation have been considered. We will pay special attention to extensions of the Tuba-Wenzl \cite{TW} representations, which exhaust all of the irreducible representations of $\B_3$ in five dimensions or less. In doing so we obtain many insights into finite-dimensional representations of $\mcLB_3$, especially in lower dimensions.

\subsection*{Results and Methodology}

 Our basic strategy is the following:

\begin{enumerate}
 \item Let $(\rho,V)$ be a representation of $\B_3$, with $\rho(\sigma_1)=A$ and $\rho(\sigma_2)=B$.
 \item Find  $S,S_1,S_2\in\End(V)$ so that $\rho(s_1)=S_1$, $\rho(s_2)=S_1S$ and $S=S_1S_2$ extends $\rho$ to a representation of $\VB_3$ (see for example Theorem \ref{existence} and Proposition \ref{BS}).
 \item Determine which pairs $(S_1,S_2)$ from the previous step factor over $\LB_3$.
\end{enumerate}

Our original goal was to carry out this procedure with an irreducible representation in step (1). Indeed, the classification results of \cite{TW} are very explicit: they show that any irreducible $\B_3$ representation of dimension $5$ or less is equivalent to a matrix representation with $A$ and $B$ in \emph{ordered triangular form}: $A$ (respectively, $B$) is upper (respectively, lower) triangular and $B_{i,i}=A_{d-i+1,d-i+1}$.  This is accomplished as part of the following theorem:
\begin{thm*}
 Let $\rho:\B_3\rightarrow \GL_d(\C)$ be an irreducible representation of $\B_3$.
 \bi
 \item If $d\leq 5$, then there exists an extension of $\rho$ to $\mcLB_3$.
 \item If $d=6$, there is an irreducible representation of $\B_3$ with no possible extensions.
 \ei
\end{thm*}
We show the above theorem using the notion of a \textbf{standard extension}, which are those for which $\rho(s_1s_2)=S=kAB$ for some $k\in\C$. The main advantage of these extensions is that they trivialize step (3) in our methodology (see Lemma \ref{mixedOK}). We have gone somewhat beyond our original goal, which we now summarize.
\subsection*{Summary of results}
\begin{enumerate}
\item
In Section \ref{General results} we characterize when standard extensions exist (Theorem \ref{existence}). We show all irreducible representations of $\mcLB_3$ when $\rho(\sigma_1)=A=B=\rho(\sigma_2)$ arise as standard extensions (Theorem \ref{AB}). We also show that extensions which factor through $\mathcal{SLB}_3$ are rare (Theorem \ref{slb}). An infinite-dimensional is given.
  \item We have determined all two-dimensional representations of $\LB_3$ in Section \ref{dim2}. In particular, we show in Theorem \ref{d2} that every irreducible two-dimensional representation of $\mcLB_3$ is a standard extension of some $\B_3$ representation and every two-dimensional $\B_3$ representation admits a standard extension.
   \item In Section \ref{dim3}, we classify extensions of irreducible three dimensional $\B_3$ representations and show that, generically, irreducible three dimensional representations of $\mcLB_3$ are standard extensions.
   \item In Section \ref{dim45}, we show all four and five-dimensional Tuba-Wenzl representations, including the reducible ones, admit standard extensions.
  \item In Section \ref{beyond}, we give an example of an irreducible $\B_3$ representation that has no possible extensions to $\mcLB_3$, and give some evidence for a conjecture that any $\B_3$ representation in the Tuba-Wenzl ordered-triangular form must have an extension.
  \end{enumerate}

\section{General results}\label{General results}
We record here our main results on extending $\B_3$ representations to
$\mcLB_3$, which are not dimension specific. We first introduce the notion of a
standard extension. Throughout this section, and for the remainder of the text,
$V$ will be a finite dimensional vector space over $\C$ unless otherwise stated.
\subsection{The standard extension}
\bd A \textbf{standard extension} of a representation $\rho:\B_3\rightarrow \GL(V)$ to $\mcLB_3$ is one for which $\rho(s_1s_2)=k\rho(\sigma_1\sigma_2)$ for some $k\in\C$.

The following lemma shows that if $A,B,S_1$ and $S_2$ satisfy relations (B1) and ($\mathfrak{S}$1) and $S_1S_2$ is proportional to $AB$, then (L1) and (L2) are also satisfied.

\ed
\bl\label{mixedOK} Let $A,B,S_1,S_2\in \End(V)$ satisfy $ABA=BAB$ and $S_1S_2S_1=S_2S_1S_2$. If $S_1S_2=kAB$ for some $k\in\C^\times$ then $ABS_1=S_2AB$ and $S_1S_2A=BS_1S_2$.
\el
\bp $ABS_1=k^{-1}S_1S_2S_1=k^{-1}S_2S_1S_2=S_2AB,\,\,\,$
$S_1S_2A=kABA=kBAB=BS_1S_2.$
\ep
Next we describe extensions of representations of the alternating group $\mathfrak{A}_3$ to $\mathfrak{S}_3$.  We use cycle notation for elements of $\mathfrak{S}_3$, and denote by $\omega$ a primitive $3$rd root of unity.

\bl\label{S3 lemma}
Let $\gamma:\mathfrak{A}_3\rightarrow \GL(V)$ be a representation given by $\gamma((1\;2\;3))=S$.  Then $\gamma$ can be extended to $\mathfrak{S}_3$ by $\gamma(s_1)=S_1$ and $\gamma(s_2)=S_2$ with $\gamma(s_1s_2)=S_1S_2=S$ if and only if $\Tr(S)\in\Z$.  Denoting by $V_\lambda$ the $\lambda$-eigenspace of $S$, we can take $S_1$ to be any involution which preserves $V_1$ and interchanges $V_{\omega}$ and $V_{\omega^{-1}}$.
 \el
\bp
If $\gamma$ is a representation of $\mathfrak{S}_3$ then $\Tr(s_1s_2)\in\Z$, since every character of $\mathfrak{S}_3$ is integral.
 Conversely, if $\Tr(S)\in\Z$ then the non-real eigenvalues $\omega^{\pm 1}$ of $S$ must appear with the same multiplicity.  Now take $S_1$ to be any involution that preserves the $1$-eigenspace of $S$ and interchanges the $\omega$ and $\omega^{-1}$ eigenspaces, and define  $S_2=S_1S$.  It is enough to verify that $S_1S_2S_1=S_2S_1S_2$, whence $S_2^2=I$ will follow from $(S_1S_2S_1)^2=S^3=I$.  We verify the equivalent condition $SS_1=S_1S^2$.  On $V_1$ this is clearly true.  For $v\in V_{\omega}$, we have $S_1S^2v=\omega^{-1}S_1v$ and $S_1v\in V_{\omega^{-1}}$ so $SS_1v=\omega^{-1}(S_1v)$.  Conversely, the relation $SS_1=S_1S^2$ with the same argument in reverse shows that the involution $S_1$ must preserve $V_1$ and interchange vectors in $V_\omega$ and $V_{\omega^{-1}}$.
 \ep

 We now characterize when a $\B_3$ representation has a standard extension.
\bt\label{trexist}\label{existence} Let $\rho:\B_3\rightarrow \GL(V)$ be a
finite-dimensional representation of $\B_3$, with $\rho(\sigma_1)=A$ and
$\rho(\sigma_2)=B$. Then the following are equivalent for $k$ and $m$ in $\C$
\begin{enumerate}
  \item[(a)] $\rho$ has a standard extension to $\mcLB_3$ with $S=kAB$ and $\Tr(S)=m$,
  \item[(b)]  $(AB)^3=k^{-3}I$ and $\Tr(AB)=k^{-1}m$ with $m\in \Z$.
   \item[(c)] $AB$ is diagonalizable and $\Tr((AB)^\ell)=k^{-\ell}m$ with $m\in\Z$ for all $\ell\leq \dim V$ not divisible by 3 and $\Tr((AB)^\ell)=k^{-\ell}\dim V$ for all $\ell\leq \dim V$ divisible by 3.

\end{enumerate}
\et
\bp The equivalence of (a) and (b) can be seen from Lemmas \ref{mixedOK} and \ref{S3 lemma}.
The implication from (a) to (c) is obtained by taking the trace of both sides of $S^\ell=(kAB)^\ell$. For (c) implies (b), we appeal to the fact that the coefficients of the characteristic polynomial of a $d\times d$ matrix $X$ is determined by the numbers $\Tr(X^\ell)$ for $1\leq\ell\leq d$ (see for example \cite{chartr1} or \cite{chartr2}).  In particular the characteristic polynomial of $AB$ is identical to the characteristic  polynomial of some matrix $X$ satisfying $X^3=k^{-3}I$ and $\Tr(X)=k^{-1}m$.  Therefore $AB$ and $X$ have the same set of eigenvalues, and since $AB$ is assumed to be diagonalizable we are done. \ep

\begin{rem}\label{explicitparam} When $\Tr(kAB)\neq0$, there is at most one possible value of $k$ for which $\Tr(kAB)\in\Z$, otherwise there are exactly three such values differing by a factor of $\omega$.
In any case, we can now produce all standard extensions of a given $\B_3$ representation satisfying either condition (b) or (c) of Theorem \ref{trexist}.
\begin{enumerate}
 \item For fixed $\rho(\sigma_1)=A$, $\rho(\sigma_2)=B$, choose $k\in\C$ so that $(AB)^3=k^{-3}I$ and $\Tr(kAB)\in\Z$.  Define $S=kAB$.
 \item Choose any $M$ so that $M^{-1}SM=(I_\ell,\omega I_t,\omega^2 I_t)$.
 \item Pick any $G\in\GL_t(\C)$ and $0\leq a\leq \ell$ and an $N\in\GL_\ell(\C)$.
 \item Set $S_1=\begin{pmatrix}
                 N(I_a,-I_{\ell-a})N^{-1} & 0 & 0\\ 0 & 0& G\\ 0 & G^{-1} &0

                \end{pmatrix}.$

 \item Define $\rho(s_1)=MS_1M^{-1}$ and $\rho(s_2)=\rho(s_1)S$.

\end{enumerate}
Notice that distinct involutions  on $\C^\ell$ with characteristic polynomial $(x-1)^a(x+1)^{\ell-a}$ are parameterized by $N(I_a,-I_{\ell-a})N^{-1}$, where $N\in\GL_\ell(\C)/(\GL_a(\C)\times \GL_{\ell-a}(\C))$.
\end{rem}

Denote by $\lceil n\rceil=\min\{k\in\Z: n\leq k\}$ the usual ceiling function.
\bprop\label{param}
 Let $\rho:\mcLB_3\rightarrow \GL(V)$ be a representation with $\rho(\sigma_1)=A$, $\rho(\sigma_2)=B$ and $\rho(s_1s_2)=S=kAB$. Suppose $W$ is invariant under $A$ and $B$. Then $W$ is $\mcLB_3$ invariant only if $\dim W\cap V_{\omega}=\dim W\cap V_{\omega^2}$. In this case, denote $\dim W\cap V_{\omega}=\dim W\cap V_{\omega^2}=m$ and let $\dim W\cap V_1=\ell$. For $\rho(s_1)=S_1$ and $\rho(s_2)=S_2$ we have:
 \begin{enumerate}
\item[(a)] $(S_1|_W,S_2|_W)$ are parametrized by
$\GL_m(\C)\times\amalg_{j=0}^{\ell}(\GL_\ell(\C)/(\GL_j(\C)\times
\GL_{\ell-j}(\C)))$.
\item[(b)] Unitary $(S_1|_W, S_2|_W)$ are parametrized by $\U_m\times
\amalg_{j=0}^{\ell}(\U_\ell/(\U_j\times \U_{\ell-j}))$.
\item[(c)] In either case, the parameter space has dimension $m^2(\lceil\frac{\ell^2-1}{2}\rceil)$ if $\ell> 1$, and $m^2$ otherwise.
  \end{enumerate}
\eprop

\begin{proof} For (a) we apply Lemma \ref{S3 lemma} to $W$ and use Remark \ref{explicitparam}. A small modification gives (b). Part (c) is a routine calculation.
\ep

\subsection{Extensions where $\rho(\sigma_1)=\rho(\sigma_2)$}
We classify here the irreducible $\LB_3$ representations $V$ in which $A=\rho(\sigma_1)=\rho(\sigma_2)=B$. In this case the mixed relation (L1) implies $S:=\rho(s_1s_2)$ commutes with $A$.  Thus the $S$-eigenspace decomposition $V=V_{1}\oplus V_{\omega}\oplus V_{\omega^2}$ must be $A$-stable.  Moreover, setting $S_1:=\rho(s_1)$ and $S_2=\rho(s_2)=S_1S$ and restricting to the subgroup $\langle s_1,s_2\rangle\cong\mathfrak{S}_3$ we see that $V_1$ and $V_\omega\oplus V_{\omega^2}$ are complementary subrepresentations of $V$.  Therefore it suffices to assume either $V=V_1$ or $V=V_{\omega}\oplus V_{\omega^2}$. We consider the first case in Theorem \ref{v1} and the second in Theorem \ref{AB}.



\bt\label{v1} Any irreducible $\LB_3$ representation $V$ with $\rho(s_1s_2)=I$ and $A=\rho(\sigma_1)=\rho(\sigma_2)=B$ and $\dim(V)>1$ is 2-dimensional of the form
$$A=\begin{pmatrix}\lambda&x\\0&-\lambda\end{pmatrix} \quad \rho(s_1)=S_1=\rho(s_2)=S_2=\begin{pmatrix}0&1\\1&0\end{pmatrix},$$
for $x\in\C$ and $\lambda\in\C^\times$. No two are isomorphic. \qed
\et

\bp
Since $V=V_1$, $S_1=S_2$.  The group $G:=\langle a,b: a^2=1,[a,b^2]=1\rangle$ is isomorphic to a split extension $K\rtimes \Z_2$ where $K$ is the kernel of $G\rightarrow \Z_2$ sending both $a$ an $b$ to $1\in\Z_2$.  Clearly $\rho(\LB_3)=\langle A,S_1\rangle$ is a quotient of $G\cong K\rtimes\Z_2$ under which the image of $K$ is $L:=\langle A^2,A S_1,S_1A\rangle$. Note that $A^2\in Z(L)$ and $AS_1=(S_1A)^{-1}$ modulo $Z(L)$. Hence we have
$1\to\Z\to L\to \Z\to1,$ which must split to give $L=\Z\times\Z$. Therefore our representation factors over $(\Z\times\Z)\rtimes\Z_2$. We show it must be at most two-dimensional and leave the rest to the reader: restrict to the abelian subgroup $\rho(K)$ and consider a one-dimensional subrepresentation $C$. By Frobenius reciprocity, $\Hom_{\mcLB_3}(\ind C,V)=\Hom_K(C, \res V)\neq0$, which means there is a non-zero map from $\ind C$ to $V$. Since $V$ is irreducible and $\ind C$ is two-dimensional then $V$ must be at most two-dimensional.
\ep

\bt\label{AB}Let $(\rho,V)$ be an irreducible representation of $\mathcal{LB}_3$ such that $A=\rho(\sigma_1)=\rho(\sigma_2)=B$ and $V_1=0$.  If $V$ is finite-dimensional, then $\dim V= 2n$ and there is a basis such that
$$
\rho(s_1s_2)=S=\mu^{-1}\omega A^2,\quad \rho(s_2)=S_2= \begin{pmatrix}0&I_n \\I_n&0\end{pmatrix},A=\diag(A_1,A_2)$$
$$A_1=\begin{pmatrix}\sqrt\mu&&&&&\\
&\begin{matrix}0&\mu \\1&0\end{matrix}&&&&\\
&&&\begin{matrix}0&\mu \\1&0\end{matrix}&&\\
&&&&\ddots&\\
&&&&&\begin{matrix}0&\mu \\1&0\end{matrix}\\
\end{pmatrix} \text{or}
 \begin{pmatrix}\sqrt\mu&&&&&\\
&\begin{matrix}0&\mu \\1&0\end{matrix}&&&&\\
&&&\ddots&&\\
&&&&\begin{matrix}0&\mu \\1&0\end{matrix}&\\
&&&&&\mp\sqrt\mu
\end{pmatrix}
$$
$$A_2=\begin{pmatrix}
\begin{matrix}0&\mu\omega \\1&0\end{matrix}&&&&\\
&&\begin{matrix}0&\mu\omega \\1&0\end{matrix}&&\\
&&&\ddots&\\
&&&&\begin{matrix}0&\mu\omega \\1&0\end{matrix}\\
&&&&&\sqrt\mu\omega^2\\
\end{pmatrix} \text{or}
 \begin{pmatrix}
\begin{matrix}0&\mu\omega \\1&0\end{matrix}&&&&\\
&\begin{matrix}0&\mu\omega \\1&0\end{matrix}&&&\\
&&&\ddots&\\
&&&&\begin{matrix}0&\mu\omega \\1&0\end{matrix}\\

\end{pmatrix}
$$\et
\bp Suppose $V=V_{\omega}\oplus V_{\omega^2}$.  Over $\C$, there is some non-zero $A^2$-eigenspace $V_{\omega,\mu}\subset V_{\omega}$. The relation $A^2S=S_2A^2S_2$ shows $S_2(V_{\omega,\mu})=V_{\omega^2,\mu\omega}$.
Both eigenspaces are preserved by $A$ and so if $V$ is irreducible, we can assume
$V=V_{\omega,\mu}\oplus V_{\omega^2,\mu\omega}$
and that there are no proper $A$ invariant subspaces of $V_{\omega,\mu}$ whose image under $S_2$ is $A$-invariant. Let $v\in V_{\omega,\mu}$ be an eigenvector of $A$ and consider the sequence $v$,  $S_2v$, $AS_2v$, $S_2AS_2v, \dots, (AS_2)^lv, S_2(AS_2)^lv,\dots$. Since $V$ is finite-dimensional, this sequence will eventually be linearly dependent. Let $n$ be the largest integer for which the first $2n$ terms are linearly independent. Note any odd number of independent terms cannot be maximal because $S_2$ is a linear bijection between those terms in $V_{\omega,\mu}$ and those in $V_{\omega^2,\mu\omega}$. Therefore the $2n+1$ term $(AS_2)^{2n}v$ is in the span of the previous $2n$ terms, which in particular means the span of the first $2n$ terms are $A$ invariant. Since it was already obviously $S_2$ and $S$ invariant, we have the span of the first $2n$ terms is all of $V$. We also have that those terms with an even number of $S_2$'s appearing span $V_{\omega,\mu}$ and those with an odd number span $V_{\omega^2,\omega\mu}$. We see with respect to this basis, the matrices have the desired form.
\ep
\begin{eg} Let $V$ be an infinite-dimensional vector space with basis $\{e_1,f_1,e_2,f_2,\dots\}$. Then we can define a representation of $\mcLB_3$ with respect to some $\mu\in\C^\times$ and $A=B$, $V_1=0$:
$$Se_i=\omega e_i, Sf_i=\omega^2f_i, S_2e_i=f_i,$$
$$Ae_1=\mu e_1, Ae_i=e_{i+1}, Af_{i-1}=f_{i}\,\text{even $i$}.$$
$$Ae_i=\mu^2e_{i-1}, Af_{i+1}=\mu^2f_{i}\,\text{odd $i$}.$$

\end{eg}

\subsection{Beyond the standard extension}
The following computational result will be useful:
\bl\label{L2} Suppose $S_1,S_2\in\GL_d(\C)$ satisfy ($\mathfrak{S}$1) and
($\mathfrak{S}$2), i.e. $S_1S_2S_1=S_2S_1S_2$ and $S_i^2=I$. Then the following
are equivalent for matrices $A,B\in\GL_d(\C)$:
\bi
\item[(a)] $ABS_1=S_2AB$,
\item[(b)] $S_2$ commutes with $ABS^{-1}$ [or equivalently with $S(AB)^{-1}$],
\item[(c)] $S_1$ commutes with $(AB)^{-1}S$ [or equivalently with $S^{-1}(AB)$],
\item[(d)] $ABS=S_2ABS_2$
\ei\el
\bp  (a) implies (b): Then $ ABS^{-1}S_2=ABS_2S_1S_2=ABS_1S_2S_1=S_2ABS_2S_1=S_2ABS^{-1}.$
(b) implies (c): $(AB)^{-1}SS_1=(AB)^{-1}S_2S_1S_2=S^{-1}S_2S(AB)^{-1}S=S_1(AB)^{-1}S$.\\
(c) implies (d): $ABS=SS^{-1}(AB)S_1S_2=SS_1S^{-1}(AB)S_2=S_2ABS_2$\\
(d) implies (a): $ABS_1=ABSS_2=S_2ABS_2S_2=S_2AB$.
\ep
\bprop\label{BS} Let $\rho$ be a $\B_3$ representation with $\rho(\sigma_1)=A$ and $\rho(\sigma_2)=B$. Suppose there is an extension $\rho^\prime(s_1)=S_1^\prime$ and $\rho^\prime(s_2)=S_2^\prime$, with $\rho^\prime(s_1s_2)=S^\prime$ to $\mcLB_3$. Suppose further that a standard extension is possible for $k\in\C^\times$. Then $\rho(s_1s_2)=S=kB^2S'$ defines an extension of $\rho$ to $\mathcal{VB}_3$. In particular, $S=k^2B^2AB$ gives $\mathcal{VB}_3$ representations.
\eprop
\bp We show $\Tr(S)\in\Z$, $SA=BS$ and $S^3=I$ and apply Lemma \ref{S3 lemma}. First we show $S$ satisfies $SA=BS$: $SA=kB^2S'A=kB^2BS'=BkB^2S'=BS.$ Now we check the trace of $S$ using Lemma \ref{L2}(d): $\Tr(S)=k\Tr(B^2S')=k\Tr(BS'A)=k\Tr(ABS')=k\Tr(S_2'ABS_2')=\Tr(kAB)\in\Z.$
Lastly, we check $S^3=I$: $S^3=k^3B^2S'B^2S'B^2S'=k^3BS'ABS'ABS'A=A^{-1}S_2'(kAB)^3S_2'A=I.$
\ep
The following well-known result helps to narrow down the possibilities for $S=\rho(s_1s_2)$, we include a proof for completeness.

\bl \label{min=char} Let $X\in\End(V)$ with $V$ a $d$-dimensional vector space.
Suppose the minimal polynomial of $X$ coincides with its characteristic polynomial. Then any $Y$ with $XY=YX$ is of the form:
$$Y=\sum\limits_{n=0}^{d-1}b_nX^n$$
\el
\begin{proof}
The hypothesis implies that there exists a $v\in V$ such that the set $\{v,Xv,\cdots,X^{d-1}v\}$ is a basis for $V$. Therefore $Yv=\sum_{n=0}^{d-1} b_nX^nv$ for some $b_n$. Since $Y$ commutes with $X^n$, we have $Y=\sum_{n=0}^{d-1}b_nX^n$ because they agree on the basis.
\end{proof}

\begin{prop}\label{uniqueness}
Let $\rho$ be any $\B_3$ representation with $\rho(\sigma_1)=A$ and $\rho(\sigma_2)=B$ such that the characteristic and minimal polynomials of $B$ coincide.  Then any extension of $\rho$ to $\mcLB_3$ has
$$\rho(s_1s_2)=S=\sum\limits_{n=0}^{d-1}a_nB^nAB$$
for some scalars $a_0, \ldots, a_{d-1}$.
\end{prop}

\begin{proof} Since $SA=BS$ by (L1) and $AB^{-1}A^{-1}=B^{-1}A^{-1}B$ by (B1), we compute:
$$(S(AB)^{-1})B=SB^{-1}A^{-1}B=SAB^{-1}A^{-1}=BSB^{-1}A^{-1}=B(S(AB)^{-1}).$$ So by Lemma \ref{min=char} $S(AB)^{-1}$  is a polynomial in $B$ of degree at most $d-1$ as required.
\ep

\subsection{Symmetric extensions}
Representations of $\LB_3$ that factor over $\SLB_3$ contain essentially no topological information, see \cite{Kane}.
We show here that it is not common to find extensions that factor through $\mathcal{SLB}_3$.
\bl\label{nilc} A matrix $S_1\in\GL_d(\C)$ commutes with a full $d\times d$
Jordan block matrix $J$ and $S_1^2=I_d$ if and only if $S_1=\pm I_d$.
\el
\bp Since $S_1$ and $J$ commute, $S_1$ is a polynomial $f(J)$ in $J=I_d+N$ where $N$ is the nilpotent Jordan matrix of nilpotency $d-1$.  Now the result follows from $S_1^2=I_d$.
\ep
\bprop\label{B2slb} Let $\rho$ be a $d$-dimensional $\mcLB_3$ representation with $\rho(\sigma_1)=A$ and $\rho(\sigma_2)=B$, whose minimal and characteristic polynomials agree, and $\rho(s_1)=S_1$ and $\rho(s_2)=S_2$. Suppose $(AB)^3=cI_d$, then $\rho$ factors through $\mathcal{SLB}_3$ if and only if $[S_2,B^2]=[S_1,A^2]=0$.
\eprop
\bp From Proposition \ref{uniqueness}, $S=a_0AB+a_1BAB+\dots+a_{d-1}B^{d-1}AB$. By Lemma \ref{L2} (b), $S_2$ commutes with $a_0+\dots+a_{d-1}B^{d-1}$. Now $B^2AB$ is a multiple of $A^{-1}B^{-1}$ when $(AB)^3=cI_d$. So we can rewrite $S(A^{-1}B^{-1})^{-1}$ up to scalar as $(a_0+\dots+a_{d-1}B^{d-1})B^{-2}$. Now $S_2$ commutes with the last expression if and only if it commutes with $B^{-2}$. Hence if we replace $A$ and $B$ by $A^{-1}$ and $B^{-1}$ in Lemma \ref{L2}, we see this is exactly (L2$^\prime$) as required. The condition for $A$ is similar.
\ep

\bt\label{slb}  Let $\rho$ be a $d$-dimensional $\B_3$ representation with $\rho(\sigma_1)=A$ and $\rho(\sigma_2)=B$, whose minimal and characteristic polynomials agree and that no eigenvalue of $A$ is the negation of another. Suppose a standard extension is possible. Then there are only finitely many extensions of $\rho$ that factor through $\mathcal{SLB}_3$.
\et
\bp By Proposition \ref{B2slb}, $S_1$ commutes with $A^2$. Since no eigenvalue of $A$ is a negation of another, $A^2$ must have a characteristic polynomial equal to its minimal polynomial. In other words, the Jordan blocks have distinct eigenvalues.  Restricting $S_1$ to a any $A^2$ eigenspace we apply Lemma \ref{nilc} to see that only finitely many $S_1$ can commute with $A^2$.  Similarly for $S_2$ and $B$.
\ep
We will see in the low-dimensional cases that the assumption on the eigenvalues cannot be removed and that finitely many cannot be strengthened to none.

\section{Two-dimensional representations}\label{dim2}
The main goal of this section will be to prove the following:

\bt\label{d2} Let $(\rho,V)$ be any two-dimensional irreducible representation
of $\mcLB_3$ with $\rho(\sigma_1)=A\neq \rho(\sigma_2)=B$, $\rho(s_1)=S_1$ and
$\rho(s_2)=S_2$. Then $\rho(s_1s_2)=-(\Tr(AB))^{-1} AB$. Moreover, every
representation of $\B_3$ with $\rho(\sigma_1)=A\neq \rho(\sigma_2)=B$ can be
extended in this way. Let $\Gr(V,1)$ be the set of one-dimensional linear subspaces of $V$.
Choices for $(S_1,S_2)$ are in one to one correspondence with $\Gr(V,1)\setminus\{V_{\omega},V_{\omega^2}\}$.
\et

\begin{rem} We only consider the $A\neq B$ case because $A=B$ was considered in the previous section. The reducible representations must satisfy $S_1=S_2$ (and therefore $A=B$) and is a straightforward calculation we have omitted. 
\end{rem}
If $A\neq B$ then we must have $S_1\neq S_2$, whence $\mathfrak{S}_3$ must act by its two-dimensional irreducible representation so that any extension to $\mathcal{LB}_3$ must also be irreducible. We also have
$$\Tr(S_1)=\Tr(S_2)=0,\Tr(S)=-1\quad\text{and}\quad\Det(S_1)=\Det(S_2)=-1.$$

\bl\label{kern} If  $\gamma((1\;2))=S_1$, $\gamma((2\;3))=S_2$ defines an irreducible two-dimensional $\mathfrak{S}_3$ representation, then the matrix equation $XS_1=S_2X$ has a two-dimensional solution space spanned by $S_1S_2$ and $S_2S_1S_2$.
\el
\bp It is easy to see that $S_1S_2$ and $S_2S_1S_2$ are solutions, which are linearly independent because $S_1\neq S_2$. Now viewing $X\rightarrow XS_1-S_2X$ as a linear transformation on $M_2(\C)$ it suffices to show the image is at least two  dimensional. Evaluating at $X=1$ and $X=S_1$ gives $S_1-S_2$ and $1-S_2S_1$ in the image. They are independent because one has zero trace and the other does not.\ep
\bl\label{trdet} If $A\neq B$ satisfy the braid relation (B1) then either
$$A=\begin{pmatrix}\lambda_1&\lambda_1\\0&\lambda_2\end{pmatrix}, B=\begin{pmatrix}\lambda_1&-\lambda_2\\0&\lambda_2\end{pmatrix}, \quad-\lambda_1/\lambda_2=\omega^{\pm1},\quad\text{or}$$
$$A=\begin{pmatrix}\lambda_1&\lambda_1\\0&\lambda_2\end{pmatrix}, B=\begin{pmatrix}\lambda_2&0\\-\lambda_2&\lambda_1\end{pmatrix},\quad\lambda_1^2-\lambda_1\lambda_2+\lambda_2^2\neq0.$$  Moreover, $(\Tr(AB))^2=\Det(AB)$.
\el
\bp
If reducible, then we are in the first case. Otherwise, we are in the second.  The last statment is a calculation.
\ep


\subsubsection*{Proof of Theorem 3.1} 
First we show $S=-(\Tr(AB))^{-1}AB$. Lemma \ref{kern} gives
$$AB=c_0S_1S_2+c_1S_2S_1S_2.$$
Taking the trace gives $c_0=-\Tr(AB)$ and taking determinant gives
$$\Det(AB)=(-\Tr(AB)+c_1)(-\Tr(AB)-c_1).$$
Hence by Lemma \ref{trdet}, $c_1=0$ and $S_1S_2=-(\Tr(AB))^{-1} AB$.

Now we show every $\B_3$ representation with $A\neq B$ admits a standard extension. By Lemma \ref{trdet}, the characteristic equation for $AB$ is
$$(AB)^2-\Tr(AB)AB+(\Tr(AB))^2=0.$$
Multiplying both sides by $AB+\Tr(AB)$, we see $(-(\Tr(AB))^{-1}AB)^3=1$.
For the last claim let $\C v$ be a one-dimensional subspace spanned by $v\in V$, then $v=v_{\omega}+v_{\omega^2}$ for unique $v_{\omega}\in V_{\omega}\setminus\{0\}$ and $v_{\omega^2}\in V_{\omega^2}\setminus\{0\}$ we define $S_1$ by $S_1(v_{\omega})=v_{\omega^2}$.

\bc\label{2slb} Suppose $V$ is an irreducible two-dimensional representation of $\mcLB_3$. Then $V$ factors through $\mathcal{SLB}_3$ if and only if $\Tr(B)=0$.
\ec
\bp Suppose $\Tr(B)=0$. Then from the characteristic equation of $B$, $B^2$ is a
scalar multiple of $I$, and hence $AB$ is a scalar multiple of $B^2AB$ which is
a scalar multiple of $A^{-1}B^{-1}$. Conversely, if the opposite relations are
satisfied, we get a standard extension of $A^{-1}$ and $B^{-1}$ and $A$ and $B$
simultaneously. So $AB$ is a scalar multiple of $A^{-1}B^{-1}$. Therefore $B^2$
is a scalar multiple of $I$ and from the characteristic equation $\Tr(B)=0$
since $B$ is not a scalar.
\ep


\section{Three dimensional representations}\label{dim3}
In this section, we show that most irreducible three dimensional representations come from standard extensions in Theorem \ref{3slb}. First we state a slightly stronger version of Theorem \ref{existence} in three dimensions.

\bl\label{3traceless} Let $\rho$ be a three dimensional representation of $\B_3$ such that $\rho(\sigma_1)=A\neq B=\rho(\sigma_2)$. Then there exists a (standard) extension to $\LB_3$ with $\rho(s_1s_2)=S=k AB$ if and only if $\Tr(AB)=\Tr((AB)^2)=0$. In this case we have $k^3=\Det(AB)^{-1}$.
\el
\bp Suppose there exists a standard extension i.e. with $S=kAB$. Since $A\neq B$ we have $S\neq I$. Restricting to $\mathfrak{S}_3$ and considering characters we see that $\Tr(S^{\pm 1})=0$. Therefore $\Tr(AB)=k^{-1}\Tr(S)=0$ and since $(AB)^2$ is a multiple of $S^{-1}$ we have $\Tr((AB)^2)=0$ as well. Conversely, suppose $\Tr(AB)=\Tr((AB)^2)=0$. The characteristic polynomial of a $3\times 3$ matrix $C$ is $$-x^3+\Tr(C)x^2+[\Tr(C)^2-\Tr(C^2)]x+\Det(C),$$ so for $C=AB$ the Cayley-Hamilton Theorem implies $(AB)^3=\Det(AB)I$.  Thus for any choice of $k^3=\Det(AB)^{-1}$ we have $(kAB)^3=I$ so that Theorem \ref{existence} implies the result.\ep

\bprop\label{SAB3} Let $\rho$ be a three dimensional representation of $\LB_3$. Suppose the minimal and characteristic polynomials of $\rho(\sigma_1)=A$ (and $\rho(\sigma_2)=B$) coincide, $\Tr(AB)=\Tr(B^2AB)=0$, $\Tr(B^4AB)\neq0$ and $A\neq B$. Then $S:=\rho(s_1s_2)=kAB$.
\eprop
\bp Observing that $\Tr(B^2AB)=\Tr((AB)^2)$, we apply Lemma \ref{3traceless} to
see that there exists a $k$ such that $(kAB)^3=I$.  Thus $(BAB)^2$ is a scalar
multiple of $I$ and so $\Tr(BAB)\neq0$ as $\rho$ is three dimensional. By Proposition \ref{uniqueness} we have
$$S=a_0AB+a_1BAB+a_2B^2AB.$$
Since $A\neq B$ and $SA=BS$, then $S\neq I$. Therefore, $\Tr(S)=0$ and taking the trace of the above expression gives $a_1=0$. Multiply by $B^2$ we obtain
$$B^2S=a_0B^2AB+a_2B^4AB.$$
Noting (by Proposition \ref{BS}) that $\Tr(B^2S)=\Tr(AB)=0=\Tr(B^2AB)$, we get $a_2\Tr(B^4AB)=0$. So $a_2=0$, i.e. $S=a_0AB$.
\ep

The following is a summary of the  results in \cite{TW} on three dimensional $\B_3$ representations.
\bl[\cite{TW}]\label{tw3} Let $\lambda_1,\lambda_2,\lambda_3\in\C^\times$ and define:
$$A=\begin{pmatrix}
  \lambda_1 & \lambda_1\lambda_3\lambda^{-1}_2+\lambda_2 & \lambda_2 \\[3mm]
  0 & \lambda_2 & \lambda_2  \\[3mm]
  0 & 0 & \lambda_3  \\
 \end{pmatrix}~~~~~~~~\text{~~~~and~~~~}~~~~~~~~
 B=\begin{pmatrix}
   \lambda_3 & 0 & 0 \\[3mm]
   -\lambda_2 & \lambda_2 & 0 \\[3mm]
   \lambda_2 & -\lambda_1\lambda_3\lambda^{-1}_2-\lambda_2 & \lambda_1
  \end{pmatrix}.$$
  Then:
  \begin{enumerate}
  \item[(a)]  $\Tr(AB)= \Tr(B^2AB)=\Tr((AB)^2)=0$ whereas $\Tr(B^4AB)=\lambda_1\lambda_2\lambda_3(\lambda_1+\lambda_2)(\lambda_1+\lambda_3)(\lambda_2+\lambda_3)$ and $(AB)^3=(\lambda_1\lambda_2\lambda_3)^2 I$.
   \item[(b)] $\rho(\sigma_1)=A$ and $\rho(\sigma_2)=B$ defines a representation of $\B_3$, which is irreducible provided $\lambda_i^2+\lambda_j\lambda_k\neq0$ for $\{i,j,k\}=\{1,2,3\}$. In this case, the minimal and characteristic polynomials of $A$ and $B$ coincide.
   \item[(c)] Let $(\rho,V)$ be any irreducible three dimensional representation of $\B_3$.  Then there exists a basis of $V$ with respect to which $\rho(\sigma_1)=A$ and $\rho(\sigma_2)=B$.
   \item[(d)] Two three-dimensional irreducible representations $\rho,\rho^\prime$ of $\B_3$ are isomorphic if and only if $Spec(\rho(\sigma_1))=\{\lambda_1,\lambda_2,\lambda_3\}=Spec(\rho^\prime(\sigma_1))$.

  \end{enumerate}
\el
We may now describe all extensions of three dimensional Tuba-Wenzl representations of $\B_3$ to $\mcLB_3$.  In particular we have a full description of extensions of \emph{irreducible} $\B_3$ representations, up to equivalence. 
\bprop\label{d3}
Let $(\rho,V)$ be a three dimensional representation of $\B_3$ such that $\rho(\sigma_1)=A$ and $\rho(\sigma_2)=B$ as in Lemma \ref{tw3}, with $Spec(A)=\{\lambda_1,\lambda_2,\lambda_3\}$.  Then $\rho$ extends to $\mcLB_3$, and either
\begin{enumerate}
 \item[(a)] $\rho(s_1s_2)=S=\gamma^{-2}AB$ where $\gamma^3=\lambda_1\lambda_2\lambda_3$, $\rho(s_1)=MS_1M^{-1}$ where $S_1=\begin{pmatrix} \pm 1&0&0\\0&0&\alpha\\
 0&\frac{1}{\alpha}&0                                                                                                                                                                                                                       \end{pmatrix}$ and $M\in\GL(V)$ is any matrix such that $\gamma^{-2}M^{-1}ABM=diag(1,\omega,\omega^2)$, or
 \item[(b)] $\rho(s_1s_2)=S\neq kAB$, $(\lambda_1+\lambda_2)(\lambda_1+\lambda_3)(\lambda_2+\lambda_3)=0$, and $\rho$ factors over $\mathcal{S}\mcLB_3$.
After permuting the eigenvalues if necessary, we have:\\
$\rho(s_1s_2)=\begin{pmatrix}
  0 & 0 & z \\[3mm]
  0 & z & z  \\[3mm]
  -\frac{1}{z^2} & \frac{1-z^3}{z^2} & -z  \\
 \end{pmatrix}$ and $\rho(s_1)=\pm\begin{pmatrix}
   1 & z-1 & z \\[3mm]
  0  & z & z  \\[3mm]
   0 & \frac{1-z^2}{z} & -z  \\
  \end{pmatrix}$ where $z\in\C^\times$ is a free parameter.
\end{enumerate}

\eprop

\bp Combining Lemma \ref{tw3} with Proposition \ref{SAB3} we get either $S=kAB$ or $(\lambda_1+\lambda_2)(\lambda_1+\lambda_3)(\lambda_2+\lambda_3)=0$.  In the first case the constant is as in Lemma \ref{3traceless} and the parametrization of $S_1$ is trivial. Hence we have (a).



Assume $(\lambda_1+\lambda_2)(\lambda_1+\lambda_3)(\lambda_2+\lambda_3)=0$ and let $S=a_0AB+a_2B^2AB$ with $a_2\neq0$ (see proof of Proposition \ref{SAB3} for why $a_1=0$). By Lemma \ref{L2}(b) $S_2=\rho(s_2)$ commutes with $S(BA)^{-1}=a_0+a_1B^2$ so $S_2$ commutes with $B^2$.  Similarly, $S_1:=\rho(s_1)$ commutes with $A^2$, so by Proposition \ref{B2slb} we see that $\rho$ factors over $\SLB_3$.  By Lemma \ref{tw3}(d), we may relabel the $\lambda_i$ so that $\lambda_3=-\lambda_2$. Set $z=\lambda_1\lambda_2(a_0+\lambda_2^2a_2)$. We obtain the required form of $S$ by solving $S^3=I$ and $\Tr(ABS)=\Tr(AB)$ and $ABS_1=S_2AB$.

\ep


One consequence of this result is that, generically, three dimensional $\B_3$ representations equi\-valent to the form of Lemma \ref{d3} do not have non-standard extensions.  This is also true in slightly greater generality, as
 the following illustrates (cf. Corollary \ref{2slb} and Theorem \ref{slb}).

\bt\label{3slb}
Let $(\rho,V)$ be an irreducible three dimensional representation of $\mcLB_3$ with $\rho(\sigma_1)=A$, $\rho(\sigma_2)=B$ with $Spec(A)=Spec(B)=\{\l_1,\l_2,\l_3\}$ and $\rho(s_1s_2)=S$. Then:
\begin{enumerate}
 \item[(a)] Any other $\LB_3$ representation $\psi$ with $\psi(\sigma_1)=A$, $\psi(\sigma_2)=B$ and $\psi(s_1s_2)=S$ is also irreducible.
 \item[(b)] Suppose the minimal and characteristic polynomials of $A$ (and $B$) coincide, $(AB)^3=cI$ and $(\lambda_1+\lambda_2)(\lambda_1+\lambda_3)(\lambda_2+\lambda_3)\neq0$. Then $\rho(s_1s_2)=kAB$.
\end{enumerate}

\et

\bp
For (a), suppose some $(\psi,U)$ satisfying the hypotheses is reducible.  Clearly $\psi(s_1s_2)=S=\rho(s_1s_2)\neq I$ since $A=SA=BS=B$ implies $\dim(V)$ is even by Theorems \ref{v1} and \ref{AB}.  Thus $S$ has 3 distinct eigenvalues $1,\omega^{\pm 1}$, and the only $\psi(s_1),\psi(s_2)$ invariant subspaces of $U$ are the $1$ eigenspace $U_1$ of $S$ and its $S$-invariant complement $U_\omega\oplus U_{\omega^{-1}}$.  So one of these two spaces is invariant under $A,B$ and $S$.  But then the same is true for $(\rho,V)$ which contradicts irreducibility as they will also be $\rho(s_1),\rho(s_2)$ invariant.

For (b), if $\B_3$ acts irreducibly on $V$, then Proposition \ref{d3} gives the result.  So assume that $\B_3$ acts reducibly on $V$ and $S=\rho(s_1s_2)\neq kAB$ for any $k$.  The hypotheses and Proposition \ref{uniqueness} imply that any $A,B$ invariant subspace $W\subset V$ is also $S$-invariant.  By Lemma \ref{3traceless} we have $S=a_0AB+a_2B^2AB$, from which it follows from Lemma \ref{L2}(b) that $S_2=\rho(s_2)$ commutes with $B^2$.  Thus if $B^2$ has distinct eigenvalues $S_2$ is a polynomial in $B^2$ by Lemma \ref{min=char} and thus $W$ is $S_2$ invariant, a contradiction.  After permuting labels we may assume $\lambda_1^2=\lambda_2^2$. We eliminate the possibility that $\lambda_1=\lambda_2$ by considering the sizes of the Jordan blocks when the minimal polynomial of $B$ must be of degree 3.  Thus $(\lambda_1+\lambda_2)=0$, a contradiction.

\ep

We illustrate our results with some applications.
\begin{eg}The Lawrence-Krammer-Bigelow representation for $\B_3$ (see for example \cite{LKB}) is defined by
$$\rho(\sigma_1)=A=\begin{pmatrix}
  tq^2 & 0 & tq(q-1) \\[3mm]
  0 &1-q & q  \\[3mm]
  0 & 1 & 0  \\
 \end{pmatrix},\quad\rho(\sigma_2)=B=\begin{pmatrix}
  1-q & 0 & 1 \\[3mm]
  0 &tq^2& tq^2(q-1)  \\[3mm]
  q & 0 & 0  \\
 \end{pmatrix},$$
for $q,t\in\C^\times$. We see that $\Tr(AB)=\Tr(B^2AB)=0$ and a standard extension is possible by Lemma \ref{3traceless}. The minimal polynomial of $A$ and $B$ both have degree 3 and the eigenvalues are not negations of each other when $tq^2\neq -1$, $tq\neq 1$ and $q\neq 1$. In this case, the irreducible extensions must be standard by Theorem \ref{3slb}.  See also \cite{bard2}, where it is shown that the Lawrence-Krammer-Bigelow representation of $\B_n$ with $n\geq 4$ does not extend except for degenerate cases.
\end{eg}
\begin{eg}
For any $1\neq t\in\C$, the $\B_3$ representation given by $\rho(\sigma_1)=A=\begin{pmatrix}
                                                    0 & t &0\\1&0&0\\0&0&1
                                                   \end{pmatrix}$ and $\rho(\sigma_2)=B=\begin{pmatrix}
                                                    1 & 0 &0\\0&0&t\\0&1&0
                                                   \end{pmatrix}$
is irreducible (see \cite{vazetal}).    Thus there is a basis with respect to which $A$ and $B$ are in the Tuba-Wenzl form.  The eigenvalues of $A$ are $\{1,\pm\sqrt{t}\}$, so that a 1-parameter family of non-standard extensions exist, by Proposition \ref{d3}(b). It is easy to verify that $\rho(s_1)=S_1=\begin{pmatrix}
                                                    0 & 1 &0\\1&0&0\\0&0&1
                                                   \end{pmatrix}$ and $\rho(s_2)=S_2=\begin{pmatrix}
                                                    1 & 0& 0\\0&0&1\\0&1&0
                                                   \end{pmatrix}$ defines an extension.
 Since $AB=\begin{pmatrix}
                                                       0 & 0&t^2\\1 &0&0\\0&1&0
                                                            \end{pmatrix}$ is is clear that $\rho(s_1s_2)\neq kAB$ for this example.  One may verify directly that this representation factors over $\mathcal{SLB}_3$.
                                                 In fact, the analogous extension works for all $n$ and factors over $\mathcal{SLB}_n$.

\end{eg}

\section{Dimensions $4$ and $5$}\label{dim45}

The following is a summary of the  results in \cite{TW} on four dimensional $\B_3$ representations.

\bl[\cite{TW}]\label{tw4} Let $\l_1,\l_2,\l_3,\l_4\in\C^\times$ and choose $\gamma^2$ such that $\gamma^4=\l_1\l_2\l_3\l_4$. Define
$$A=\begin{pmatrix}
  \lambda_1 & (1+\dfrac{\lambda_1\lambda_4}{\gamma^2}+\dfrac{\lambda_1^2\lambda_4^2}{\gamma^4})\lambda_2 & (1+\dfrac{\lambda_1\lambda_4}{\gamma^2}+\dfrac{\lambda_1^2\lambda_4^2}{\gamma^4})\lambda_3 & \lambda_4 \\[5mm]
  0 & \lambda_2 & (1+\dfrac{\lambda_1\lambda_4}{\gamma^2})\lambda_3 & \lambda_4  \\[5mm]
  0 & 0 & \lambda_3 & \lambda_4  \\[3mm]
  0 & 0 & 0 & \lambda_4
 \end{pmatrix}$$
 and
 $$B=\begin{pmatrix}
   \lambda_4 & 0 & 0 & 0 \\[3mm]
   -\lambda_3 & \lambda_3 & 0 & 0 \\[3mm]
   \dfrac{\lambda_2^2\lambda_3}{\gamma^2} & -(\dfrac{\lambda_2\lambda_3}{\gamma^2}+1)\lambda_2 & \lambda_2 & 0 \\[5mm]
   -\dfrac{\lambda_1\lambda_2^3\lambda_3^3}{\gamma^6} & (\dfrac{\lambda_2^3\lambda_3^3}{\gamma^6}+\dfrac{\lambda_2^2\lambda_3^2}{\gamma^4}+\dfrac{\lambda_2\lambda_3}{\gamma^2})\lambda_1 & -(\dfrac{\lambda_2^2\lambda_3^2}{\gamma^4}+\dfrac{\lambda_2\lambda_3}{\gamma^2}+1)\lambda_1 & \lambda_1
  \end{pmatrix}$$
\begin{enumerate}
\item Setting $\rho(\sigma_1)=A$ and $\rho(\sigma_2)=B$ defines a representation of $\B_3$, and every irreducible $\B_3$ representation of dimension $4$ is equivalent to such a representation.

\item $\Tr(AB)=-\gamma^2$ and $(AB)^3=-\gamma^{6}I_4$.
\end{enumerate}
\el
\bprop All four dimensional Tuba-Wenzl representations admit a standard extension. In particular, all irreducible four dimensional representations of $\B_3$ admit a standard extension.
\eprop
\bp By Lemma \ref{tw4}, $(-\gamma^{-2}AB)^3=I_4$ and $\Tr(-\gamma^{-2}AB)=1$. Now apply Theorem \ref{existence}.
\ep
\bprop Any extension of a four dimensional Tuba-Wenzl representation satisfies $\Tr(S)=1$.
\eprop
\bp $\Tr(S)$ must either be $-2$ or $1$ since $S\neq1$.
If $\Tr(S)=-2$, then $e_1+Se_1+S^2e_1=0$ and $e_2+Se_2+S^2e_2=0$ means $S$ must have diagonal entries $0,0,-1,-1$ and other entries zero unless they are on the skew diagonal, in which case they would be $x,x,-x^{-1},-x^{-1}$. We see that such a matrix cannot satisfy $SA=BS$. Equating the $(2,4)$ entry gives $\lambda_4=-\lambda_3$. Then equating the $(3,4)$ entry gives $x(x-1)=\lambda_2^2/\gamma^2$.
\ep

\begin{eg} If $V$ and $W$ are both two-dimensional irreducible representations of $\mcLB_3$ then $V\otimes W$ is a four dimensional representation with $S=kAB$.
\end{eg}

The following is a summary of the classification of simple $\B_3$ representations of dimension 5 found in \cite{TW}.
\bl[\cite{TW}]\label{tw5}Let $\l_1,\l_2,\l_3,\l_4,\l_5\in\C^\times$ with $\gamma$ a fixed fifth root of $\l_1\l_2\l_3\l_4\l_5$. Define

$$A=\begin{pmatrix}
  \lambda_1 & (1+\dfrac{\gamma^2}{\lambda_2\lambda_4})(\lambda_2+\dfrac{\gamma^3}{\lambda_3\lambda_4}) & (1+\dfrac{\lambda_1\lambda_5}{\gamma^2})(\lambda_3+\gamma+\dfrac{\gamma^2}{\lambda_3}) & (1+\dfrac{\lambda_2\lambda_4}{\gamma^2})(\lambda_3+\dfrac{\gamma^3}{\lambda_2\lambda_4}) & \dfrac{\gamma^3}{\lambda_1\lambda_5} \\[5mm]
  0 & \lambda_2 & \lambda_3+\gamma+\dfrac{\gamma^2}{\lambda_3} & \lambda_3+\gamma+\dfrac{\gamma^3}{\lambda_1\lambda_5} & \dfrac{\gamma^3}{\lambda_1\lambda_5} \\[5mm]
  0 & 0 & \lambda_3 & \lambda_3+\dfrac{\gamma^3}{\lambda_1\lambda_5} & \dfrac{\gamma^3}{\lambda_1\lambda_5} \\[5mm]
  0 & 0 & 0 & \lambda_4 & \lambda_4 \\[2mm]
  0 & 0 & 0 & 0 & \lambda_5 \\
 \end{pmatrix}$$
 and $B$ by $B_{i,j}=(-1)^{i-j}A_{6-i,6-j}.$
\begin{enumerate}
\item Defining $\rho(\sigma_1)=A$ and $\rho(\sigma_2)=B$ defines a representation of $\B_3$, and every irreducible $\B_3$ representation of dimension $5$ is equivalent to such a representation.

\item $(\gamma^{-2}AB)^3=I_5$ and $\Tr(\gamma^{-2}AB)=-1$.
\end{enumerate}

\el
From this lemma and Theorem \ref{existence} we immediately have:
\bprop All five-dimensional Tuba-Wenzl representations admit a standard extension. In particular, all irreducible five-dimensional representations of $\B_3$ admit a standard extension.
\eprop

\bprop Any extension of a five-dimensional Tuba-Wenzl representation satisfies $\Tr(S)=-1$.
\eprop
\bp Assume otherwise, $\Tr(S)=2$. In this case the non real eigenspace of $S$ is a unique two-dimensional subspace. So the span of $e_1, Se_1,S^2e_1$ and $e_2,Se_2,S^2e_2$ must have a two-dimensional intersection. Assuming the skew triangular form of $S$, this determines the 3rd row and column completely (all but (3,3) entry zero, which is 1). Applying $SA=BS$ to this row and column gives two of the skew diagonal entries next to the (3,3) entry are both -1. Now $S^3e_1=e_1$ and $S^3e_2=e_2$ gives the matrix for $S$ has only non-zero entries on the skew diagonal $x,-1,1,-1,x^{-1}$ and the last row and column: $x, v, 0, v, 1$, $x^{-1},-v^{-1},0,-v^{-1},1$. This now should be checked to fail the relation $SA=BS$.
\ep

In dimension $d=4$ and $5$, every irreducible $\B_3$ representation can be extended to a standard $\mathcal{LB}_3$ representation. Now we show that most extensions are standard.
That is, whenever an extension exists, $S=kAB$ for some $k\in\C^\times$ unless the eigenvalues $(\lambda_1,\ldots,\lambda_d)$ of $A$ and $B$ are zeros of a set of polynomials.

\bprop\label{triangle lemma} Let $\rho$ be an extension to $\LB_3$ of an irreducible $\B_3$ representation of dimension $4$ or $5$ in the form of Lemma \ref{tw4} or \ref{tw5}, with $\rho(\sigma_1)=A$, $\rho(\sigma_2)=B$ and $\rho(s_1s_2)=S$.  Then $AB$, $S$ and $BSA$ are all skew lower triangular.  Furthermore, $S^2$ and $(BSA)^2$ are skew upper triangular matrices.
\eprop
\bp A direct computation shows that $AB$ is skew lower triangular. Since the characteristic and minimal polynomials of $A$ coincide, Proposition \ref{uniqueness} implies $S$ is the product of a lower triangular matrix by a skew lower triangular matrix, which is evidently skew lower triangular. Therefore, $S^2=S^{-1}$ is skew upper and the same holds for $BSA=B^2A$ by Proposition \ref{BS}.
\ep

\bt\label{unique}  Let $\rho$ be  an irreducible $d$-dimensional $\B_3$-representation with $d=4$ or $5$ as in Lemmas \ref{tw4} and \ref{tw5} with $\rho(\sigma_1)=A$ and $\rho(\sigma_2)=B$.
  For each $d=4$ and $5$, there is a variety $V(J_d)$ where $J_d\subset K[\lambda_1,\ldots,\lambda_d,\gamma]$, such that if $(\lambda_1,\ldots,\lambda_d)\in \mathbb{C}^d\setminus V(J_d)$, then $$S=kAB$$
  is the only solution for the loop braid relations, where $k=-\gamma^{-2}$ for $d=4$ and $k_5=\gamma^{-2}$ for $d=5$.
\et

\begin{proof} By Proposition \ref{uniqueness}, $S=\sum_{i=0}^{d-1}b_iB^iAB$. We will prove that, except when $(\lambda_1,\ldots,\lambda_d,\gamma)$ is a solution to a certain system of polynomial equations, all coefficients $b_i$ vanishes except $b_0$.
Let $P_{ij}$ be the $(i,j)$-entry of $S^2$ and $Q_{ij}$ be the $(i,j)$-entry of $(BSA)^2$. Then $P_{ij}$ and $Q_{ij}$ are second degree homogeneous polynomials of $b_0,\ldots,b_{d-1}$. Proposition \ref{triangle lemma} gives
$$P_{ij}=0 \text{~~~~and~~~~} Q_{ij}=0 \text{,~~~~if~~~~}i<j \eqno (*)$$
Note that for $i<j$, the term $b_0^2$ does not appear in $P_{ij}$ and $Q_{ij}$ for $(AB)^2$ and $(BABA)^2$ are skew upper triangular. Therefore, for $i<j$, System ($*$) can be viewed as a homogeneous linear system of unknowns $b_mb_n$. It has more equations than unknowns when $d\geq 2$. In fact, it has $d^2-d$ equations and $N_d$ unknowns $b_mb_n$ ($m+n>0$), where $N_d=(d+2)(d-1)/2$ and $N_d\leq d^2-d$ for $d \geq 2$. In general, it only has zero solutions depending on whether the coefficient matrix has rank $N_d$.

Let $M_d$ be the coefficient matrix of System ($*$) and $J_d$ be the set of the determinants of all $N_d\times N_d$ submatrices of $M_d$. If $(\lambda_1,\ldots,\lambda_d,\gamma)\notin V(J_d)$, then System ($*$) only has zero solutions, that is $$b_mb_n=0\text{~~~~for~~~~}m+n>0.$$

Therefore, for generic values of $(\lambda_1,\ldots,\lambda_d,\gamma)$ (i.e. away from $V(J_d)$ we must have $S=b_0AB$.\end{proof}

\section{Beyond dimension $5$}\label{beyond}
We show that six-dimensional irreducible representations of $\B_3$ do not always extend.

\bprop Let $\omega=e^{2\pi i/3}$ be a primitive 3rd root of unity and define:
$$A=\begin{pmatrix}
  1 & -\omega+1 & -\omega^2+1 & \omega-1 & \omega^2-1 & \omega-1  \\[3mm]
  \omega^2-1 & \omega^2 & 0 & -\omega^2+1 & 0 & 0  \\[3mm]
  \omega^2-1 & \omega^2-1 & \omega^2 & -\omega^2+1 & -\omega^2+1 & 0  \\[3mm]
  0 & \omega-1 & \omega^2-1 & -\omega & -\omega^2+1 & -\omega+1 \\[3mm]
  -\omega^2+1 & -\omega^2+1 & 0 & \omega^2-1 & -1 & 0  \\[3mm]
  -\omega^2+1 & -\omega^2+1 & -\omega^2+1 & \omega^2-1 & \omega^2-1 & -1
 \end{pmatrix}$$
 and
 $$B=\begin{pmatrix}
   1 & -\omega+1 & -\omega^2+1 & -\omega+1 & -\omega^2+1 & -\omega+1  \\[3mm]
   \omega^2-1 & \omega^2 & 0 & \omega^2-1 & 0 & 0  \\[3mm]
   \omega^2-1 & \omega^2-1 & \omega^2 & \omega^2-1 & \omega^2-1 & 0  \\[3mm]
   0 & -\omega+1 & -\omega^2+1 & -\omega & -\omega^2+1 & -\omega+1 \\[3mm]
   \omega^2-1 & \omega^2-1 & 0 & \omega^2-1 & -1 & 0  \\[3mm]
   \omega^2-1 & \omega^2-1 & \omega^2-1 & \omega^2-1 & \omega^2-1 & -1
  \end{pmatrix}.$$
  Then $\rho(\sigma_1)=A$ and $\rho(\sigma_2)=B$ define an irreducible representation of $\B_3$ which cannot be extended to an $\mathcal{LB}_3$ representation.
  \eprop
  \bp
 These two matrices follows the construction in \cite{L}, where the representation is shown to be indecomposable.  We verify (using {\texttt Magma} \cite{MAGMA}) that the dimension of the algebra ge\-nerated by $A,B$ is 36, we see that the representation is irreducible.  Alternatively, one may use \cite[Remark 2.11(4)]{TW} to verify irreducibility: in this case $\Det(A)^{6}=1$ so that if there were a non-trivial $r<6$-dimensional subrepresentation $W\subset \C^6$, the eigenvalues of $A$ on $W$ must satisfy $(\mu_1\cdots\mu_r)^{36}=1$.  Computing the characteristic polynomial one finds that the eigenvalues of $A$ are $e^{\pi i/3}$ and the $5$ roots of an irreducible $5$ degree polynomial in $\Q(\sqrt{3}i)[x]$.  In particular, the minimal polynomial of $B$ (and $A$) coincides with the characteristic polynomial. In any case, any subrepresentation either has dimension $1$ or a $1$ direct complement. One then verifies that $A$ and $B$ do not have a common eigenvector of eigenvalue $e^{\pi i/3}$.

Now if an extension of $\rho$ to $\mathcal{LB}_3$ exists, then, by Proposition \ref{uniqueness}, $S=\sum_{i=0}^5b_iB^iAB$.  The requirement $S^3=I$ and a calculation imply that $S=qAB$ or $S=qB^2AB$ where $q$ is any third root of unity. But for such $S$, $tr(S)\notin\mathbb{R}$. Therefore, this irreducible $\mathcal{B}_3$ representation cannot be extended to a $\mathcal{LB}_3$ representation.
\ep

In \cite[Subsection 1.5]{TW} there is a general construction of $(d+1)\times (d+1)$ representations in ordered triangular form.  We recall the details for the reader's convenience.  For $0\leq i\leq d$ define $\overline{i}=d-i$ and let $c\in\C^\times$. For $0\leq i\leq d$, let $\lambda_i$ satisfy $\lambda_i\lambda_{d-i}=c$.  Define matrices $A$ and $B$ by
$A_{ij}=\begin{pmatrix} \overline{i}\\ \overline{j}\end{pmatrix}\lambda_{\overline{i}}$ and $B_{ij}=(-1)^{i+j}\begin{pmatrix} {i}\\ {j}\end{pmatrix}\lambda_i$.  The identity:
$\sum_{k=0}^d(-1)^{k}\begin{pmatrix} {i}\\ {k}\end{pmatrix}\begin{pmatrix} \overline{k}\\ \overline{j}\end{pmatrix}=\begin{pmatrix} {i}\\ {j}\end{pmatrix}$ is useful to verify the relation $ABA=BAB$ so that these matrices define a $\B_3$ representation.  One computes:
$(AB)_{ij}=c(-1)^{\overline{j}}\begin{pmatrix} {i}\\ \overline{j}\end{pmatrix}$ and $(AB)^3=(-1)^{d}c^3I$.  Define $S=\frac{(-1)^{d}}{c}AB$, that is: $$S_{ij}=(-1)^{d+\overline{j}}\begin{pmatrix} {i}\\ \overline{j}\end{pmatrix}=(-1)^{j}\begin{pmatrix} {i}\\ \overline{j}\end{pmatrix}.$$ Clearly $\Tr(S)\in\R$ (in fact $\Tr(S)\in\{0,\pm 1\}$ and $S^3=I$) and this representation admits a standard extension.
Using \cite[Remark 2.11(4)]{TW} (due to Deligne) it is possible to show that, for sufficiently generic eigenvalues, these $\B_3$ representations are irreducible.

Our investigations of low-dimensional $\mcLB_3$ representations suggest the following:
\begin{conj}
 Suppose that $\rho$ is an irreducible $d$-dimensional matrix representation of $\B_3$ such that
 $\rho(\sigma_1)=A$ and $\rho(\sigma_2)=B$ are in ordered triangular form (\cite{TW}), that is, $A$ is upper triangular and $B$ is lower triangular with $B_{i,i}=A_{d-i+1,d-i+1}$.
Then $\rho$ has a standard extension to $\mcLB_3$.
\end{conj}
To prove this conjecture it is enough to show that for some root $k$ of $cx^3-1$ where $(AB)^3=cI$ we have $\Tr(kAB)\in\Z$.  Fixing some root $c^{1/3}$, we have $Spec(c^{1/3}AB)\subset\{1,\omega^{\pm 1}\}$ where $\omega=e^{2\pi i/3}$.  If these appear with multiplicities $\mu_1,\mu_{\pm}$ then there is a choice of $k$ such that $\Tr(kAB)\in\Z$ if and only if two of these eigenvalues coincide.

\end{document}